\documentclass [11pt]{article}
\usepackage{amsmath, amsthm, amssymb}
\usepackage{graphicx}
\newfont{\bbb} {msbm10}
\newcommand{\R}{\Bbb{R}}
\newcommand{\bS}{\Bbb{S}}

\newcommand{\Z}{\Bbb{Z}}

\newcommand{\D}{\Bbb{D}}

\newcommand{\sbs}{\subset}
\newcommand{\ra}{\rightarrow}

\newcommand{\bq}{{\bar{q}}}

\newcommand{\dsigma}{{\dot{\sigma}}}

\newcommand{\p}{\partial}

\newcommand{\cD}{{\cal{D}}}
\newcommand{\cS}{{\cal{S}}}

\newcommand{\cC}{{\cal{C}}}

\newcommand{\cA}{{\cal{A}}}

\newcommand{\cB}{{\cal{B}}}

\newcommand{\B}{\Bbb{B}}

\newcommand{\s}[1]{{\sf{#1}}}

\newcommand{\0}[1]{_{_{#1}}}
\newcommand{\rC}{{\rm{C}}\,}

\newcommand{\sL}{{\sf{Link}}}

\usepackage[all]{xy}

\begin{document}

\title{Normal Smoothings for Smooth Cube Manifolds}
\author{Pedro Ontaneda\thanks{The author was
partially supported by a NSF grant.}}
\date{}

\maketitle

\begin{abstract} A smooth cube manifold $M^n$ is a smooth $n$-manifold $M$ together with a smooth cubulation on $M$. (A smooth cubulation is similar to a smooth triangulation, but with cubes instead of simplices).
The cube structure provides, for each open
$k$-subcube $\dsigma^k$, rays that are normal to $\dsigma$. Using these rays we can
construct {\it normal charts} of the form $\D^{n-k}\times\dsigma\ra M$, where we are identifying
$\D^{n-k}$ with the cone over the link of $\dsigma$ (these identifications are
arbitrary and the identification between $\p \D^{n-k}$ and the link of $\dsigma$ is called a {\it link smoothing}). These normal charts respect the product
structure of $\D^{n-k}\times\dsigma$ and the radial structure of $\D^{n-k}$.
A complete set of normal charts gives a (topological) {\it normal atlas} on $M$. If this
atlas is smooth it is called a {\it normal smooth atlas} on $M$ and induces a {\it
normal smooth structure} on $M$ (normal with respect to the cube structure).
In this paper we prove that every smooth cube manifold has a normal smooth
structure, which is diffeomorphic to the original one.
This result also holds for smooth all-right-spherical manifolds.
\end{abstract}

\noindent {\bf \large  Section 0. Introduction.}

As mentioned in the Abstract, in this paper we address the following
quite natural question. Given a a cube complex (very popular objects these
days)
one can construct ``normal charts" in an obvious natural way (charts respect
radial directions). Question: given a smooth cube complex (smooth manifold
+compatible cube structure), does there exist a normal smooth structure?
(i.e is there a smooth atlas formed of normal charts?) In this paper we
give an affirmative answer to this question.\\

The results in this paper are key ingredients to smooth the metric of a strictly hyperbolized manifold (see
\cite{O2} and \cite{O1}). 
This smoothing implies the existence of a large class of negatively curved
Riemannian manifolds, showing that, in some sense, Riemanian negative curvature abounds in nature.\\

 Here is a brief explanation how 
normal smooth structures on cube manifolds are used in \cite{O2}.
To any manifold $M$ with a smooth cubulation $K$ one can associate a smooth piecewise hyperbolic manifold $K_X$ produced via the strict hyperbolization procedure of Charney and Davis (based on Gromov's hyperbolization
process). Using a smooth normal atlas $\cA$ on $M$ with respect to
$K$, we can contruct a similar ``normal" smooth atlas $\cA_{K_X}$ on the strictly hyperbolized manifold $K_X$ (this is done in detail in \cite{O1}, and also 
briefly explained in \cite{O2}).  Because of the form of normal charts (they preserve some of the ``cube structure"),
it is possible to use the atlas $\cA_{K_X}$  to smooth the natural piecewise hyperobolic metric
on $K_X$ via local warped product constructions. \\

We next give a more detailed description of normal structures and state the Main Theorem.\\

For the basic definitions and results about cube and spherical  complexes see
for instance \cite{BH}. Recall that a spherical
complex is an {\it all-right-spherical complex} if all of its edge lengths are equal to 
$\pi/2$. 
Given a (cube or all-right-spherical) complex $K$ we use the same notation $K$ for the complex itself (the collection of all
closed cubes or simplices) and its realization (the union of all cubes or simplices).
For $\sigma\in K$ we denote its interior by $\dsigma$.
In this paper we assume that all cube or all-right-spherical complexes satisfy the following condition: any two closed cubes or
simplices intersect
on at most one (possibly empty) common subcube or subsimplex.
\vspace{.1in}

Let $M^n$ be a smooth manifold of dimension $n$. Recall that a {\it cubulation} of  $M$ is given by $(K,f)$, where 
$K$ is a cube complex and $f:K\ra M$ is a homeomorphism.
The cubulation is {\it smooth} if
$f$  a non-degenerate
$PD$ homeomorphism \cite{MunkresLectures}, that is, for all $\sigma\in K$  we have $f|_{\sigma}$ is a smooth embedding. Sometimes we will write $K$ instead
of $(K,f)$. The smooth manifold $M$ together with a smooth cubulation is
a {\it smooth cube manifold}. A {\it smooth all-right-spherical manifold} is
defined analogously.\vspace{.1in}

In this paper
we consider $\sL(\sigma^{j},K)$,  the (geometric) link of an open $j$-cube or $j$-all-right simplex $\sigma^j$,
as the union of the end points of straight segments of small length $\epsilon>0$
emanating perpendicularly (to $\dsigma^j$) from some fixed point  $x\in \dsigma^j$. 
We say that the link is {\it based at $x$.} And the (open) star $\s{Star}(\sigma,K)$ as the union of such half-open $\epsilon$-segments based at x. We can identify the star with the (open) cone of the link
$\rC\sL (\sigma,K)$  (or $\epsilon$-cone) defined as
$\rC \, \sL(\sigma, K)=\sL(\sigma,K)\times [0,\epsilon)\,/\, \sL(\sigma,K)\times\{0\}$.
 We shall denote the {\it cone point} by $o$ or, more specifically, by 
$o\0{\rC \, \sL(\sigma, K)}$.
Thus a point $x$ in $\rC \, \sL(\sigma, K)$, different from the cone point $o$, can be written as
$x=t\,u$, $t\in (0,\epsilon)$, $u\in \sL(\sigma, K)$. For $s>0$ we get the
{\it cone homothety} $x\mapsto sx=(st)u$ (partially defined if $s>1$). 
If we want to make explicit the dependence of the link or the cone on $\epsilon$ we shall write $\sL_\epsilon(\sigma,K)$ or $\rC _\epsilon\,\sL(\sigma,K)$ respectively. 
Also, we will always take $\epsilon <1/2$ ($<\pi/4$ in the spherical case) and it can be verified
that all results in this paper (unless otherwise stated) are independent
of the choice of the $\epsilon$'s. 
As usual we shall identify the $\epsilon$-normal neighborhood
of $\dsigma$ in $K$ with $\rC _\epsilon\,\sL(\sigma,K)\times
\dsigma$ (or just $\rC \,\sL(\sigma,K)\times
\dsigma$). Hence a cone homothety induces a {\it neighborhood homothety}
obtained by crossing it with the identity $1\0{\dsigma}$.
\vspace{.1in}

 Recall that the link $\sL(\sigma^i,K)$,  $\sigma^i\in K$, has a natural
all-right piecewise spherical structure, which induces a simplicial structure and thus a $PL$ structure on  $\sL(\sigma^i,K)$.\vspace{.1in}

From now on we assume that $K$ is a $PL$ manifold,
and $f:K\ra L$ a homoemorphism.
(For instance, $K$ is a smooth cubulation, or triangulation,
of $M$; in this case the $PL$ structure on $M$ induced by $K$ is Whitehead compatible with $M$, hence $K$ is a $PL$ manifold, see \cite{KiSi}, p.10.) Therefore 
the link $\sL(\sigma^i,K)$ is $PL$ homeomorphic to $\bS^{n-i-1}$,
for every $\sigma^i\in K$.
A {\it link smoothing for}  $\dsigma^i$ (or $\sigma^i$)
is just a homeomorphism $h_{\sigma^i}:\bS^{n-i-1}\ra\sL(\sigma^i,K)$.
The {\it cone} of $h_{\sigma^i}$ is the map 
$$\rC \,h_{\sigma^i}:\D^{n-i}\longrightarrow \rC\sL (\sigma^i, K)$$
\noindent given by $t\,x= [x,t]\mapsto t\,h_{q^i}(x)=[h_{q^i}(x),  \,t]$,
 where we are canonically identifying the $\epsilon$-cone  of $\bS^{n-i-1}$ with the open disc $\D^{n-i}$. We remark that we are not assuming $h_{\sigma^i}$ to be $PL$.\vspace{.1in}

A link smoothing $h_{\sigma^i}$  induces the following  smoothing of the normal neighborhood of
$\dsigma^i$:
$$h^\bullet_{\sigma^i}=f\,\, \circ\,\,\Big(\rC \,h_{\sigma^i}\times 1_{\dsigma^i}\Big):\D^{n-i}\times \dsigma^i\longrightarrow M $$
The pair $(\,h^\bullet_{\sigma^i}\, ,\,\D^{n-i}\times \dsigma^i\,)$, or simply 
$h^\bullet_{\sigma^i}$, is a {\it normal
chart} on $M$. Note that the collection  $\cA=\big\{\,(\,h^\bullet_{\sigma^i}\, ,\,\D^{n-i}\times 
\dsigma^i\,)\,\big\}_{\sigma^i\in K}$ is a topological atlas for $M$.
Sometimes will just write $\cA=\big\{\,h^\bullet_{\sigma^i}\,\big\}_{\sigma^i\in K}$.
The topological atlas $\cA$ depends uniquely on the
the complex $K$, the map $f$ and the collection of link smoothings 
$\{h_{\sigma}\}_{\sigma\in K}$. 
To express the dependence of the atlas on the set of links smoothings
we shall write $\cA=\cA\big(\{h_{\sigma}\}_{\sigma\in K}\,\big)$
(this is different from $\cA=\big\{\,h^\bullet_{\sigma^i}\,\big\}_{\sigma^i\in K}$,
as written above).\vspace{.1in}

The most important feature about these normal
atlases is that they preserve the radial and sphere (link) structure given by $K$.
These features make normal atlases very powerful tools for geometric
constructions.\vspace{.1in}

Note that not every collection of link smoothings induce a smooth atlas, but when 
the atlas is smooth we call $\cA$ a {\it normal smooth atlas on $M$ with respect to} $K$ and the corresponding smooth structure $\cS'$
a {\it normal smooth structure on $M$ with respect to $K$}. 
Notice that in this case the maps $f|_{\dsigma^i}:\dsigma^i\ra (M,\cS')$ are 
smooth embeddings.
Also note that $\cA$ is smooth if and only if there is a smooth structure
$\cS'$ such that all normal charts 
$h^\bullet_{\sigma^i}:\D^{n-i}\times \dsigma^i\longrightarrow (M,\cS') $
are smooth embeddings. Here is our main result.\vspace{.1in}

\noindent {\bf Main Theorem.} {\it Let $M$ be a smooth cube manifold,
with smooth structure $\cS$.
Then $M$ admits a normal smooth structure $\cS'$ diffeomorphic to $\cS$.}\vspace{.1in}

Specifically,   if $M^n$ is a smooth manifold with smooth structure $\cS$ and $K$ is a smooth cubulation of $M$, then we can choose link smoothings $ h_{\sigma^i}$, for all $\sigma^i\in K$, such that the atlas $\cA=\cA\big(\{h_{\sigma}\}_{\sigma\in K}\,\big)$
is smooth. The smooth atlas $\cA$ is normal with respect to $K$. Moreover the normal smooth structure $\cS'$ (with respect to $K$),
induced by $\cA$, is diffeomorphic to $\cS$.\vspace{.1in}

\noindent {\bf Addendum to Main Theorem.} {\it The statement of the Main
Theorem also holds for smooth all-right-spherical complexes.}
\vspace{.1in}

\noindent {\bf Remarks.}

\noindent {\bf 1.} 
It can be checked from the proof of the Main Theorem that in fact $\cS'$ is isotopic to $\cS$ for 
$n\neq 4$. That is, there is a diffeomorphism $\phi : (M,\cS)\ra (M,\cS')$ that is (topologically) isotopic to the identity map $1_M$. With some care we can probably include the case $n=4$.

\noindent {\bf 2.} Note that the image of the chart $h^\bullet_{\sigma}$
is the open normal neighborhood $\stackrel{\circ}{\s{N}}_\epsilon(\dsigma,K)$
of width $\epsilon$ of $\dsigma$ in $K$.
Even though we are assuming, for simplicity, that $\epsilon<1/2$ ($\epsilon<\pi/4$
in the spherical case) it can be checked from the proof of the Main Theorem that we can
actually take $\epsilon=1$ ($\epsilon=\pi/2$) for the charts. 

\noindent {\bf 3.} The proof of the Main Theorem 
uses an inductive construction of a normal atlas.
The idea is to keep adding normal charts. The main ingredient
in the construction is Lemma 1.2 which says that the links
of cubes are smooth submanifolds. Lemma 1.2 does not work
if we replace cubes by (linear) simplices, or even by equilateral
simplices: a key ingredient in the proof of Lemma 1.2. is the smoothness of the function in
formula (1), but in the equilateral simplicial case
the analogous formula has a function which is not smooth. Of course, if we ``reparametrize" the rays in the
linear case then the proof works, because this is just the
all-right-spherical case.

\noindent {\bf 4.} Every smooth cubulation induces a
smooth triangulation (by subdivision) and every  smooth triangulation induces a smooth cubulation (by canonically ``cubifying" the
standard simplices). Therefore every smooth manifold admits
a smooth cubulation, hence by the Main Theorem, every
smooth manifold admits a normal smooth structure, relative
to some cube structure.

\noindent {\bf 5.} As application of Lemma 1.2 we get a simple
 proof of the following fact: the vanishing of
$\Theta_n$ (the group of homotopy $n$-spheres), $n\leq 6$,
$n\neq 4$, implies that every $PL$ $n$-manifold is smoothable,
for $n\leq 7$. (See remark before Proposition 1.4.)

\noindent {\bf 6.} The Main Theorem is of course not true
if we replace ``smooth cubulation" by just ``cubulation".
In particular the Main Theorem is not true if $K$ is not
a $PL$ manifold.

\noindent {\bf 7.} Let $\cS'$ be as in the Main Theorem. As mentioned above the maps 
$f|_{\dsigma^i}:\dsigma^i\ra (M,\cS')$ are smooth embeddings (recall $\dsigma$ is an {\sf open} simplex), but
we can not expect the maps $f|_{\sigma^i}$ to be embeddings
(recall $\sigma$ is a {\sf closed} simplex). Nor we can expect the map $f:K\ra (M,\cS')$ to be PD. That is, we cannot expect
$(K,f)$ to be a smooth cubulation of $(M,\cS')$. This does not even happen in the next simple example.
(On the other hand, since $\cS$ and $\cS'$ are diffeomorphic
$(K,\psi)$ is a smooth cubulation of $(M,\cS')$ for some other map
$\psi$. In fact this map $\psi$ can also be constructed rather explicitly by induction using Lemma A.2.2 in the appendix; see A.2.5. Lemma A.2.2 corrects the lack
of differentiability showed in the next example.)\vspace{.1in}

\noindent{\bf Example.} Consider $\R^2$ with its canonical cube structure, that is, the vertices of the cubes are points in $\R^2$
with integer coordinates. Thus the 2-cubes  are translations
of  $[0,1]^2$.
Note that the canonical smooth structure of $\R^2$ is
normal with respect to $K$: choose the normal charts
$h^\bullet_{\sigma^1}$ of 1-cubes
to be the obvious product charts, and the link smoothings $h_{\sigma^0}$
for the vertices  to be translations of the inclusion $\bS^1=\bS^1_\epsilon\ra \R^2$.
Consider now $\R^2$ with a cube structure given by pulling back the canonical one using
the map $p\0{\ell}:\R^2\ra\R^2$,  $p\0{\ell}(r,\theta)= (r,\ell \,\theta)$, where we are using polar coordinates, and $\ell\neq \pm 1$ is an integer. 
Note that $p\0{\ell}$  is not differentiable at the origin.
Let $K_\ell$ be this new cube complex and $f:K_\ell\ra \R^2$ be the inclusion. Hence the images of the open cubes $f(\dsigma')$, $\sigma'\in K_\ell$ are the components of $p\0{\ell}^{-1}(\dsigma)$, $\sigma\in K$.
Notice that the canonical smooth structure of $\R^2$
is also normal  with respect to $K_\ell$: (1) choose the normal charts of an open cube $\dsigma$ in $K_\ell$, $\sigma\neq 0$ ( i.e. a cube different from the origin) to be the composition of $p\0{\ell}^{-1}$ with the the normal chart 
of $p\0{\ell}(\sigma)$, that is, choose $h_{\sigma}^\bullet=p\0{\ell}^{-1}\circ h_{p\0{\ell}(\sigma)}^\bullet$ (this is well defined), (2) for the origin choose the link smoothing $h\0{\ell}$ of the link $\sL(0,K\0{\ell})$ of the origin as the lifting by
$p\0{\ell}$ of the link smoothing $h$ of $\sL(0,K)$ (which is just the inclusion):
$$\begin{array}{ccc}\bS^1& \stackrel{h\0{\ell}}{\ra}&\sL(0,K\0{\ell})
\\ \downarrow &&\downarrow p\0{\ell}\\
\bS^1& \stackrel{h}{\ra}&\sL(0,K)\end{array}
$$
\noindent where the left vertical arrow is the map $z\mapsto
z^l$, $z\in\bS^1$.
Then for any closed 2-cube $\sigma^2\in K_\ell$ containing the origin we have that  $f|_{\sigma^2}:\sigma^2\ra \R^2$ 
coincides with $p\0{\ell}|_{\sigma^2}$ which is not differentiable at 
the origin, hence $f$ is not $PD$.\vspace{.1in}

So, even though the normal charts $h^\bullet_{\sigma^i}$ have potent geometric
properties,  the example above shows that the $PL$ structure given by $K$ has to be sacrificed. (But we can replace it by
an equivalent one).\vspace{.1in}

In section 1 we prove the Main Theorem and in section 2 we generalize our results to the case of manifolds with
discrete point singularities. There is one appendix.\\

We are grateful to the referee for his/her comments and suggestions.
\vspace{.3in}

\noindent {\bf Section 1. Proof of Main Theorem.}

To simplify our notation we will write the proof for cubical complexes $K$; the proof for the all-right-spherical complex case 
is similar. We will denote an open $i$-cube by $q^i$ and the corresponding closed cube
by $\bq^i$. We write $\sL(q,K)=\sL(\bq,K)$.
Also we denote by $K_j$ the $j$-skeleton of $K$, i.e.  the union
of all cubes of dimension $\leq j$. Also write $M_j=f(K_j)$. We will prove the Theorem by induction on $k=n-j$. Consider the following
statement:

\begin{enumerate}
\item[{\bf S($k$)}]
{\it There is a smooth structure $\cS_{k}$ on $M-M_j$ and link smoothings $h_{q^i}$ for all $i$-cubes, $i> j $, such that
the corresponding normal charts $h^\bullet_{q^i}:\D^{n-i}\times 
q^i\longrightarrow \big(M-M_j,\cS_{k}\big)$ are smooth embeddings.}
\end{enumerate}

\noindent {\bf Remark.} If {\bf S($k$)} holds, then $\big\{\,(\,h^\bullet_{q^i}\, ,\,\D^{n-i}\times 
q^i\,)\,\big\}_{\bq^i\in K,\, i>j}$ is a smooth atlas for $(M-M_j,\cS_k)$.
Moreover, the smooth structure $\cS_{k}|_{M-M_{j+1}}$ coincides with $\cS_{k-1}$, provided the link smoothings of both structures coincide for all the $q^i$, $i>j+1$.
Also  the maps $f|_{q^i}:q^i\ra (M-M_j,\cS_{k})$ are embeddings, $i>k$.

The following is a key ingredient in the proof of the Main Theorem
and is the reason why everything works.\vspace{.1in}

\noindent {\bf Lemma 1.2.} {\it Assume {\bf S($k$)} holds. Then  $f\big(\sL(q^{j},K))$ is a smooth submanifold of $(M-M_j,\cS_{k})$,
for every $j$-cube $q^{j}$ of $K$.}\vspace{.1in}

In general $f\big(\sL(q^{j},K))$ may not be diffeomorphic to
$\bS^{k-1}$, but if it is, we have the following addition to Lemma 1.2.\vspace{.1in}

\noindent {\bf Addendum to Lemma 1.2.} {\it Let $h_{q^j}:\bS^{k-1}\ra \sL(q^j,K)$ be a link smoothing such that $f\circ h_{q^j}$
is a diffeomorphism. Then the corresponding normal chart
$$h^\bullet_{q^j}|_{(\D^k-0)\times q^j}:(\D^k-0)\times q^j\ra \big(M-M_j,\cS_k\big)$$ 
\noindent is a smooth embedding.}\vspace{.1in}

\noindent {\bf Proof of Lemma 1.2.} By the remark above it is enough to verify that $f(\sL(q^{j},K))$ is a smooth submanifold in
every chart  $h^\bullet_{q^i}:\D^{n-i}\times 
q^i\longrightarrow \big(M-M_j,\cS_{k}\big)$, that is, we have to show that $(h^\bullet_{q^i})^{-1}(f(\sL(q^{j},K)))$ is a submanifold
of $\D^{n-i}\times  q^i$. We can assume that $\bq^j$ is a subcube of $\bq^i$.
To be specific consider $\sL(q^{j},K)=\sL\0{\epsilon}(q^{j},K)$,
$\sL(q^{i},K)=\sL\0{\delta}(q^{i},K)$ and
assume that $\sL(q^{j},K)$ is based at $x\in q^j$.  Write $S=(h^\bullet_{q^i})^{-1}
\big(f(\sL(q^{j},K))\big)$. By the definition of $h^\bullet_{q^i}$
we have $S=(\rC h^\bullet_{q^i}\times 1_{q^i})^{-1}\big(\sL(q^{j},K)\big)$. Also write
 $\bq^i=\bq^l\times \bq^j$. Then $S\sbs \D^{n-i}\times q^l\times\{x\}$. Identify $\bq^l$ with
$\bq^l\times\{ x\}$ and $x$ with a vertex of $\bq^l$.
By the definition of (geometric) link, we have that $S$ is the set of points  $(a,b)\in\D^{n-i}\times q^l$ 
such that the segment $[x,(\rC h_{q^i}(a),b)]$ has length $\epsilon$ (and is perpendicular to $q^j$) in $K$. 
Since $\rC h_{q^i}\times 1_{q^i}$ is, by definition, radial in the first coordinate and the identity on the second, $S$ is the set of points $(a,b)\in\D^{n-i}\times q^l$ that satisfy the equation 

\begin{equation*}\epsilon^2\,=\,\delta^2\,d_{\D^{n-i}}^2(a,0)\,+\,d_{\bq^l}^2(b,x)\tag{1}\end{equation*}

\noindent where $d_{\D^{n-i}}$ and $d_{\bq^l}$ are the euclidean distances on $\D^{n-i}\sbs\R^{n-i}$ and
$\bq^l\sbs\R^l$, respectively. 
Therefore $S$ is a submanifold of 
$\D^{n-i}\times q^l$. This proves the lemma.\vspace{.1in}

Now we prove the addendum. 
Recall that we are identifying a neighborhood of $q^i$ in $K$
(or $f(q^i)$ in $M$) with $\rC \sL(q^{i},K)\times q^i$
(respectively $f\Big(\rC \sL(q^{i},K)\times q^i\Big)$).

Now we claim that the homothety map $$f\Big(\Big(\rC \sL(q^{j},K)-o\Big)\times q^j\Big)\times (0,1)\ra f\Big(\Big(\rC \sL(q^{j},K)-o\Big)\times q^j\Big),\,\,\,\,\,\,\,f(x,y,s)\mapsto f(sx,y)$$
\noindent is smooth, where we are considering 
$f\Big(\Big(\rC \sL(q^{j},K)-o\Big)\times q^j\Big)$ with the structure
$\cS_k$. To see that the homothety above is smooth recall that we are writing $\bq^i=\bq^l\times \bq^j$ and notice that the above map
written in the chart $h^\bullet_{q^i}$ is $(a, b,c,s)\mapsto
(sa, sb  ,c)$, $(a,b,c,s)\in \R^{n-i}\times\R^{l}\times\R^j\times (0,1)$. (Here we are identifying $\bq^i=\bq^l\times \bq^j$ with $[0,1]^i
=[0,1]^l\times [0,1]^j$ so that $q^j$ corresponds to
$\{0\}\times [0,1]^j$.) This shows that the homothety above 
smooth. \vspace{.1in}

We also need the fact that the distance function $f(x,y)\mapsto |x|$, is smooth. Here $(x,y)\in (\rC \sL(q^{j},K)-o)\times q^j$,
and $|x|$ is the distance to the vertex in $\rC \sL(q^{j},K)$. To see this let $q^i$ so that $(x,y)\in q^i$. Using the notation in the previous paragraph we have the square of the distance function in the
$h^\bullet_{q^i}$ chart is $(a,b,c)\mapsto \delta^2\,d_{\D^{n-i}}^2(a,0)\,+\,d_{\bq^l}^2(b,x)$ (as in equation (1)). This shows the square
of the distance function (hence the distance function itself because
it does not vanish) is smooth.

To prove the addendum just note that for $x\neq 0\in \D^k$ we have \vspace{.05in}

\hspace{1.85in}$h^\bullet_{q^j}(x,y)\,=\,f\,\Big(\,|x|\,h_{q^j}\big(\frac{1}{|x|}x    \big)\,,\,
y   \,\Big)$\vspace{.05in}

\noindent where we are using the homothety and the distance function in the right side
of the equality.
Hence $h^\bullet_{q^j}$ is smooth (outside $q^j$). Again using
charts it is straightforward to show that $h^\bullet_{q^j}$ is an immersion (outside $q^j$).
Since $h^\bullet_{q^j}$
is a (topological) embedding it follows that it is a smooth embedding. This proves the addendum.\vspace{.1in}

\noindent {\bf Remarks.} 

\noindent {\bf 1.} Equation (1) in the all-right spherical
case can be obtained using the spherical law of cosines. The equation is
$$
\epsilon\,=\, \cos\big( \delta\,d_{\D^{n-i}}(p,0)\big)
\,\cos\big(\,d_{\bq^l}(q,x) \big).
$$
\noindent {\bf 2.} Even if $f\big(\sL(q^{j},K))$ is diffeomorphic to
an exotic sphere $\Sigma^{k-1}$ the Addendum to Lemma 1.2 remains true with one change: consider now $\D^k-0$ with the smooth structure given by the identification
$\D^k-0=\Sigma^{k-1}\times(0,1)$.\vspace{.1in} 

Lemma 1.2 implies that the smooth structure $\cS_k$ induces a smooth structure on the sphere $f\big(\sL(q^{j},K)\big)$, but the
smooth structure and the $PL$ structure (induced by $K$)  are not necessarily compatible (see remark 7 in the introduction). If $k\leq 7$, $k\neq 5$
the $(k-1)$-sphere $f\big(\sL(q^{j},K)\big)$,   with the induced smooth structure, is diffeomorphic to the canonical sphere (see \cite{MK}). The next lemma says
that the same is true for $k=5$.\vspace{.1in}

\noindent {\bf Lemma 1.3.} {\it Let $j=n-5$ and consider the 4-sphere $f\big(\sL(q^{j},K)\big)$ with the induced smooth structure
from Lemma 1.2. Then  $f\big(\sL(q^{j},K)\big)$ is diffeomorphic to $\bS^4$.}\vspace{.1in}

The proof of this lemma is given in the Appendix, but here is
a quick sketch of the proof. Recall that $f\big(\sL(q^{n-5},K)\big)$
 with the $PL$ structure induced by $K$ is a $PL$ 4-sphere
 (because $K$ induces a $PL$ triangulation on $M$). But,
 as mentioned before, the $PL$ structure on $f\big(\sL(q^{n-5},K)\big)$ induced by $\cS_5$ may not coincide with the $PL$
 structure induced by $K$. The trick here is to prove that
 these two $PL$ structures are $PL$ equivalent. The way
 to prove this is similar to the proof (given below, see
 condition {\bf C}$(k)$) that all structures $\cS_k$ are
 diffeomorphic to the original $\cS$. \vspace{.1in}

We now prove the Main Theorem.  Recall we are writing $j+k=n$.
We begin with $k=1$ and define $\cS_1=\cS|_{M-M_{n-1}}$.
Then  $\cA_1=\big\{\,(\,h^\bullet_{q^n}\, ,\,\D^{0}\times 
q^n\,)\,\big\}_{\bq^n\in K}$ is a smooth atlas for $\cS_1$.
Using lemmas 1.2 and 1.3 (for the case $k$=5) we can construct $\cS_k$, $k\leq 7$ inductively:
$\cS_{k+1}$ has an atlas $\cA_{k+1}$ which  is obtained from $\cA_k$ by adding the
charts  $\big\{\,(\,h^\bullet_{q^{j}}\, ,\,\D^{n-j}\times 
q^j\,)\,\big\}_{\bq^j\in K}$, for some smoothings $h_{q^j}$ of the links of the $j$-cubes,
such that $h_{q^j}:\bS^{k-1}\ra f(\sL(q^j,K))$ is a diffeomorphism
(which exists by 1.2 or 1.3).
Hence, by construction, we have that $\cA_k= \big\{\,(\,h^\bullet_{q^i}\, ,\,\D^{n-i}\times 
q^i\,)\,\big\}_{\bq^i\in K,\, i>j}$, for $k\leq 7$.\vspace{.1in}

\noindent {\bf Remark.} As mentioned in the introduction
 every smooth manifold admits
a smooth cubulation. Therefore Lemmas 1.2, 1.3
and the simple construction of the atlas $\cA_k$ given above
prove the following fact: the vanishing of
$\Theta_n$ (the group of homotopy $n$-spheres), $n\leq 6$,
$n\neq 4$, implies that every $PL$ $n$-manifold is smoothable,
for $n\leq 7$. Of course if the Smooth Poincare Conjecture
in dimension 4 was true then we would not need Lemma 1.3.\vspace{.1in}

The construction of the atlas $\cA_k$ given above proves the Main Theorem for $n\leq 3$, since (diffeomorphism classes of) smooth structures are unique for $n\leq 3$. For $n=4$ it only remains to prove that the smooth
structures $\cS'=\cS_5$ and $\cS$ are diffeomorphic, and this is done in the appendix 
(see A.2.3).
Hence from now on we assume $n\geq 5$. We will prove by induction the following stronger statement. Let
{\bf S'($k$)}, $k\geq 4$, be the statement obtained from {\bf S($k$)} by adding the extra condition {\bf C($k$)} where: \\

\noindent {\bf C($4$):}\,\,\, {\it $\cS_{4}$ extends to a smooth structure $\cS_{4}'$ on the whole $M$, 
and $\cS_{4}'$ is diffeomorphic to $\cS$. }\\

\noindent {\bf C($k>4$):}
\,\,\, {\it  $\cS_{k}$ extends to a smooth structure $\cS_{k}'$ on the whole $M$, 
and $\cS_{k}'$ is diffeomorphic to $\cS_{k-1}'$.}\\

Recall that we have proved {\bf S($k$)} for $k\leq 7$. The next proposition
shows our strategy to prove the Main Theorem.\vspace{.1in}

\noindent {\bf Proposition 1.4.} {\it We have that}
\begin{enumerate}
\item[{\bf (1)}] {\it Statement} {\bf C}(4) {\it holds}. 
\item[{\bf (2)}] {\it Statements} {\bf C}($k$) {\it hold for $k=5,\, 6$}.
\item[{\bf (3)}] {\it Statement} {\bf S'}($k-1$) {\it implies} {\bf S}($k$).
\item[{\bf (4)}] {\it Statements} {\bf C}($k$) and {\bf S}($k+1$)
{\it imply} {\bf C}($k+1$) {\it for $k\geq 6$.}
\end{enumerate}
 \noindent {\bf Proof of Main Theorem assuming Proposition 1.4.}\\
Just note that it follows from (2), (3), (4) and the fact that we have already proved     {\bf S($7$)} that 
statement {\bf S'}($k$) implies {\bf S'}($k+1$), for $k\geq 6$. Thus induction
shows that statement {\bf S'}($n+1$) is true. This proves the Main Theorem 
assuming Proposition 1.4.\vspace{.1in}

\noindent {\bf Proof of Proposition 1.4.}

First we need a definition. A {\it spindle neighborhood map} of
an open cube $q^i$ in $K$ is a topological embedding $\alpha: \D^{n-i}\times q^i\ra K$ such that:  (1) $\alpha(0,x)=x$ and \, (2) the diameters of
the fibers $\alpha (\{x\}\times\D^{n-i})$ (with respect to the path metric on $K$) tend uniformly to zero as $x$
tends to $\p q^i$. (This terminology is essentially due to J. Munkres \cite{Munkres}). 
The image of $\alpha$ is a {\it spindle neighborhood} of $q^i$. 
\vspace{.1in}

The proof of  (1) uses some ingredients of the proof of Lemma 1.2, and
it is given in the second part of the appendix.\vspace{.1in}

We prove (2) next. Assume {\bf C}($k-1$) holds, where $k=5,\, 6$.
(Here we use that (1) holds.) Choose a spindle neighborhood for each $(j+1)$-cube $q^{j+1}$. We can assume all spindle neighborhoods
to be disjoint.
Since we are constructing $\cS_k$ and $\cA_k$ by induction, the smooth structure $\cS_{k}|_{M-M_{j+1}}$ coincides with $\cS_{k-1}$.
Hence the smooth structures $\cS_k$ and $\cS_{k-1}'$ coincide on $M-M_{j+1}$.
Let $q^{j+1}$ be a $(j+1)$-cube and pull back these two smooth structures to $ \D^{k-1}\times
q^{j+1}$ via the corresponding
spindle map $\alpha_{q^{j+1}}$. Call these structures $\cB$ and $\cC$.
They coincide outside $q^{j+1}=\{0\}\times
q^{j+1}$. \\

Now, the Classification Theorem 10.1 of \cite{KiSi}, p.194,
says that concordance classes of smooth structures
on $ \D^{k-1}\times q^{j+1}$, modulo the complement $C$
of an open neighborhood of $q^{j+1}$ in $\D^{k-1}\times
q^{j+1}$, is in one-to-one correspondence with
$$\big[\,\big(\D^{k-1}\times q^{j+1}\,,\,C\,\big)\,\,,\,\,TOP/O\big]\,=\,
\pi_{k-1}(TOP/O)$$
From now on we chose $C$ to be  a small spindle
neighborhood of $q^{j+1}$ in $ \D^{k-1}\times q^{j+1}$.
Since $\pi_{k-1}(TOP/O)=\pi_{k-1}(K(\Z_2,3))=0$, 
for $4\leq k-1 \leq 6$ (\cite{KiSi} p.201), we have that $\cB$ and $\cC$
are concordant modulo $C$. Therefore (see Theorem 4.1 on p.25 of \cite{KiSi}) $\cB$ and $\cC$ are isotopic modulo $C$. Consequently there is a diffeomorphism 
$\beta:( \D^{k-1}\times q^{j+1}, \cB)\ra ( \D^{k-1}\times q^{j+1},\cC)$, which is the identity on $C$. 
Moreover $\beta$ is isotopic to the identity modulo $C$.
Using these maps $\beta$ we can define a diffeomorphism
$(M-M_j,\cS_k)\ra (M-M_j,\cS_{k-1}')$ that extends to a homeomorphism $\phi : M\ra M$
by declaring $\phi(x)=x$, for $x\in M_j$. Take now $\cS_k'$ to be the pull back of $\cS_{k-1}'$ by $\phi$.
This proves (2).\vspace{.1in}

To prove (3) we need the following result. Recall we are assuming
$n\geq 5$.\vspace{.1in}

\noindent {\bf Lemma 1.5.} {\it Assume {\bf S'($k$)} holds.
Consider the $(k-1)$-sphere $f(\sL(q^{j},K))$ with the induced smooth structure
from Lemma 1.2. Then  $f(\sL(q^{j},K))$ is diffeomorphic to $\bS^{k-1}$.}\vspace{.1in}

\noindent {\bf Proof.} 
 Let $h_{q^j}:\bS^{k-1}\ra \sL(q^j,K)$
be a homeomorphism. 
Pull back the smooth structure  $\cS'_k$ by the map $h^\bullet_{q^j}$ to
obtain a smooth structure $\cB$ on $\D^{n-j}\times q^j$. 
The Addendum to Lemma 1.2 (also see remark 2 before Lemma 1.3) implies that the smooth structure $\cB$ is a product structure outside 
$q^j=\{ 0\}\times q^j\sbs \D^k\times q^j$. By Lemma 1.3 we can assume 
$k-1> 6$ and we can apply the Product Structure Theorem (see Theorem 5.1 on p.31 of \cite{KiSi}) to obtain a diffeomorphism
$\big( \D^k\times q^j, \cC\times \cS_{\R^j} \big)\ra\big(\D^k\times q^j,\cB \big)$
for some smooth structure $\cC$ on the disc $\D^k$. 
Here $\cS_{\R^j}$ is the canonical smooth structure on $\R^j$. Furthermore, we can assume this diffeomorphism to be the identity outside a small neighborhood of $q^j$.
Hence the smooth structure $\cB|_{\bS^{k-1}}=\cC|_{\bS^{k-1}}$ extends to a smooth structure
on the whole disc. We can apply now the smooth $h$-cobordism Theorem to conclude that
$\cB|_{\bS^{k-1}}$ is diffeomorphic to the canonical structure on the sphere. This
proves the lemma.\vspace{.1in}

We can prove (3) now:
just use Lemma 1.5 and the addendum to Lemma 1.2 to construct $\cA_k$ and $\cS_k$ as done before. 
Therefore {\bf S'($k-1$)} implies {\bf S($k$)}. \vspace{.1in}

It remains to prove (4), that is, we have to prove that  {\bf C($k-1$)} and {\bf S($k$)}
imply  {\bf C($k$)}, i.e. that $\cS_k$ extends to the whole $M$. 
The proof is similar to the proof of (2), but first we
may have to change the smoothings of the links of $j$-cubes (given in the proof of Lemma 1.5).
The following lemma will indicate how to change the smoothing $h_{q^j}$.
Recall that for $f:\bS^k\ra\bS^k$ the map $\rC f:\D^{k+1}\ra\D^{k+1}$
is the cone of $f$ defined by $\rC f (x)=|x|f(\frac{x}{|x|})$, $x\neq 0$ and
$f(0)=0$.\vspace{.1in}

\noindent {\bf Lemma 1.6.} {\it Let $\cD$ be a smooth structure on the disc $\D^k$, $k \geq 6$. Assume $\cD$ coincides with the
canonical smooth structure outside a small neighborhood
of the origin.
Then there is a diffeomorphism $g:\bS^{k-1}\ra\bS^{k-1}$ such that the smooth structure
$(\rC \,g)^*\cD$ is diffeomorphic to the canonical one, by a diffeomorphism that is the identity outside a small neighborhood
of the origin.}\vspace{.1in}

\noindent {\bf Proof.} Let $G:\D^k\ra(\D^k,\cD)$ be a diffeomorphism. Since, by hypothesis,  $\cD$ coincides with the canonical smooth structure outside a small neighborhood
of the origin, we can assume that $G$ is radial outside
a small neighborhood of the origin (to see this use the uniqueness of collars or deformation of vector fields techniques). Take $g=G|_{\bS^{k-1}}$
and note that $((\rC g)^{-1}\circ G)^*(\rC g)^*\cD=G^*\cD$ is the canonical smooth
structure on $\D^k$.  Moreover $(\rC g)^{-1}\circ G$ is the identity outside
a small neighborhood of the origin.
This proves Lemma 1.6.\vspace{.1in}

We now prove (4). Assume {\bf C($k$)} and {\bf S($k+1$)} hold.
Let $\cB=\cB(h_{q^j})$ be the pull-back smooth structure $(h^\bullet_{q^j})^*\cS'_{k}$. 
On the other hand note that, by construction, the pull-back smooth structure $(h^\bullet_{q^j})^*\cS_{k+1}$ is the canonical 
smooth structure on $\D^k\times q^j$. Since $\cS_k'$ extends
$\cS_k$, and $\cS_k$ and $\cS_{k+1}$ coincide on $M-M_j$, we have that $\cB$ is a smooth structure on $\D^k\times q^j$ that coincides with the canonical one
outside $q^j$.

For a self-homeomorphism $g$ on $\bS^{k-1}$ 
write $g^\bullet=(\rC \, g)\times 1_{q^j}$, which is a self-homeomorphism 
on $\D^{k-1}\times q^j$.\vspace{.1in}

\noindent {\bf Claim.} {\it There is a self-homeomorphism $g$ on $\bS^{k-1}$
such that $\cB( h_{q^j}\circ g)=(g^\bullet)^*\cB$ is diffeomorphic to the canonical
smooth structure on $\D^k\times q^j$, via a diffeomorphism that is the
identity outside a small spindle neighborhood of $q^j$ in $\D^k\times q^j$.}\vspace{.1in}

\noindent {\bf Proof of claim.} Recall that we are denoting by $\cS_{\R^j}$ the canonical smooth structure on $\R^j$.
Using the Product
Structure Theorem (see Theorem 5.1 in \cite{KiSi}, p.31), we get a diffeomorphism
$$\gamma: \big( \D^k\times q^j ,\cD\times\cS_{\R^j}\big)\ra \big( \D^k\times q^j,\cB \big)$$
\noindent for some smooth structure $\cD$ on the disc $\D^k$. Furthermore,
$\gamma$ is the identity outside a small spindle neighborhood of $q^j$ and isotopic to the identity (by an isotopy that is constant outside a small spindle neighborhood of $q^j$). Apply Lemma 1.6 to $\cD$ to obtain a diffeomorphism $g$. 
Note that $(g^\bullet)^{-1}\circ\gamma\circ g^\bullet$ is a diffeomorphism between the smooth structures $(\gamma\circ g^\bullet)^*\cB$ and $(g^\bullet)^*\cB$, and this diffeomorphism is the identity on the complement of a small neighborhood of $q^j$ .
But by Lemma 1.6 we have that $(\gamma\circ g^\bullet)^*\cB=(g^\bullet)^*(\cD\times \cS_{\R^j})=(\rC g)^*\cD\times\cS_{\R^j}$
which is diffeomorphic to the canonical smooth structure on $\D^k\times q^j$,
via a diffeomorphism that is the identity on the complement of a small neighborhood of $q^j$. 
Therefore we have proved that 
$(g^\bullet)^*\cB$ is diffeomorphic to the canonical
smooth structure on $\D^k\times q^j$, via a diffeomorphism that is the identity outside a small neighborhood of $q^j$ in 
$\D^k\times q^j$ (not necessarily spindle). Call this diffeomorphism
$\eta$, that is, $\eta^*(g^\bullet)^*\cB$ is diffeomorphic to the canonical
smooth structure, and $\eta$ is the identity outside a small neighborhood of $q^j$. We have to change $\eta$ by one that
is the identity outside a small {\sl spindle} neighborhood of $q^j$.

To be specific let $\eta$ be the identity outside $B\times q^j$,
where $B\sbs\D^k$ is an open ball centered at 0 of small radius.
Let $\theta$ be a self diffeomorphism of
$\D^k\times q^j$ and $V$ a small spindle neighborhood
of $q^j$ such that: (1) $\theta$ is of the form $\theta (u,v)=(\theta_v(u), v)$, that
is, $\theta$ maps the fibers $\D^k\times\{ v\}$ to themselves,
(2) $V=\coprod_{v\in q^j} B_v\times \{ v\}$, where $B_v\sbs \D^k$
is an open ball centered at 0 of small radius $r\0{v}$,
and $r\0{v}\ra0$ uniformly as $v\ra\p q^j$,
(3)  $\theta$ is the identity inside a very small spindle neighborhood
$U$ of $q^j$, with $\bar{U}\sbs V$, (4) $B\sbs\theta_v (B_v\times \{ v\})$. By (3) and the fact that $\cB$ coincides with the canonical structure outside $q^j$ we have that $\theta$ is also a self 
diffeomorphism of $( \D^k\times q^j,\cB \big)$. It is now
straightforward to verify that  $\theta^{-1}\circ\eta\circ\theta$ is a
diffeomorphism between
$(g^\bullet)^*\cB$ and the canonical
smooth structure on $\D^k\times q^j$, and that this map is the identity outside the small spindle neighborhood $V$ of $q^j$. This proves the claim.\vspace{.1in}

From the claim it follows that
if we replace the link smoothing $h_{q^j}$ by  the new link smoothing $ h_{q^j}\circ g$ and take $\cB=\cB(h_{q^j}\circ g)$
instead of $\cB=\cB(h_{q^j})$
we get that the new structure $\cB$ is now diffeomorphic to the canonical one
 modulo the complement of a small spindle neighborhood of $q^j$. That is
there is a diffeomorphism $\beta=\beta_{q^j}:\D^k\times q^j\ra (\D^k\times q^j,\cB)$ that
is the identity outside a small spindle neighborhood of $q^j$.
We can now 
proceed as in the the proof of (2) (i.e. the case $k\leq 7$) and use the maps $\beta_{q^j}$ to obtain
a diffeomorphism  
$(M-M_{j-1},\cS_{k+1})\ra (M-M_{j-1},\cS_{k}')$ that extends to a homeomorphism $\phi : M\ra M$
by declaring $\phi(x)=x$, for $x\in M_j$. Take now $\cS_{k+1}'$ to be the pull-back of $\cS_{k}'$ by $\phi$. This proves (4) and completes the proof of Proposition 1.4
\vspace{.3in}

\noindent {\bf 2. Manifolds with codimension zero singularities.}

In this section we treat the case of manifolds with a one point singularity.
The case of manifolds with many (isolated) point singularities
is similar.\vspace{.1in}

Let $Q$ be a smooth manifold with a one point singularity $q$, that is,
$Q-\{q\}$ is a smooth manifold and there is a topological embedding
$\rC_1 N\ra Q$, with $o\0{\rC N}\mapsto q$, which is a smooth embedding outside 
the vertex $o\0{\rC N}$. Here $N=(N,\cS_N)$ is a closed smooth manifold (with smooth
structure $\cS_N$).
Also $\rC_1 N $ is the closed cone of $N\times [0,1]/N\times\{0\}$. 
We write $\rC_1N\sbs Q$ and sometimes
we shall identify
$\rC_1 N-\{o\0{\rC N}\}$ with $N\times (0,1]$. 
We say that the {\it singularity $q$ of $Q$ is  modeled on  $\rC N$.}
Also we say that $(K,f )$ is a {\it smooth cubulation of $Q$} if

\begin{enumerate}
\item[(i)] 
$K$ is a cubical complex.
\item[(ii)] $f:K\ra Q$ is a homeomorphism. Write $f(p)=q$ and $L=\sL(p,K)$.
\item[(iii)] $f|_{\sigma}$ is a smooth embedding for every
cube $\sigma$ not containing $p$.
\item[(iv)]  $f|_{\sigma-\{p\}}$ is a smooth embedding for every
cube $\sigma$ containing $p$.
\item[(v)]  
$L$ is $PL$ homeomorphic to $(N,\cS_N)$.
\end{enumerate}

Many of the definitions and results given in previous sections still hold (with minor changes) in the case of manifolds with a one point singularity:
\begin{enumerate}
\item[{\bf (1)}] A {\it link smoothing} for $L=\sL(p,K)$ (or $p$) is just a
homeomorphism $h_p:N\ra L$.
A {\it set of  link smoothings for $K$} is a set of link smoothings for the sphere links plus
a link smoothing for $L$.
\item[{\bf (2)}] Given a set of link smoothings for $K$ we get a set of normal charts
as before. For the vertex $p$ we mean the cone map 
$h_p^\bullet=f\circ\rC h_p:\rC N\ra Q$. We will also denote the restriction
of $h_p^\bullet$ to $\rC N-\{o\0{\rC N}\}$ by the same notation $h_p^\bullet$.
As before $\{h^\bullet_\sigma\}_{\sigma\in K}$ is a {\it (topological) normal atlas on $Q$
with respect to $K$}. The atlas on $Q$ is {\it smooth} if all transition functions are smooth, where for the case $h_p^\bullet:\rC N-\{o\0{\rC N}\}\ra Q-\{q\}$ we are identifying $\rC N-\{
o\0{\rC N}\}$ with $N\times (0,1]$ with the product smooth structure obtained from
\s{some} smooth structure ${\tilde{\cS}}_N$ on $N$.
A smooth normal atlas on $Q$ with respect to $K$ induces, by restriction,  a smooth normal structure on $Q-\{q\}$ with respect to $K-\{p\}$ (this makes sense even though
$K-\{p\}$ is not, strictly speaking, a cube complex).
\item[{\bf (3)}] We say that the set $\{h_\sigma\}$ is {\it smooth} 
if the atlas $\cA=\{h^\bullet_\sigma\}_{\sigma\in K}$ is smooth.  In this case we say that
the smooth atlas  $\cA$ (or the induced smooth structure, or the set $\{h_\sigma\}$) 
is {\it correct with respect to $N$} if $\cS_N$ and ${\tilde{\cS}}_N$ are diffeomorphic.
\end{enumerate}

The Main Theorem  also holds in this context:\\

\noindent {\bf Theorem 4.1.} {\it  Let $Q$ be a smooth manifold
with one point singularity $q$ modeled on $\rC N$, where $N$ is a closed
smooth manifold. Let $(K,f)$ be a smooth cubulation of $Q$. Then
$Q$ admits a normal smooth structure with respect to $K$, which restricted to $Q-\{q\}$ is diffeomorphic to $Q-\{q\}$.
Moreover this normal smooth structure is correct with respect to $N$ if}
\begin{enumerate}
\item[(a)] {\it $dim\, N\leq 4$.}
\item[(b)] {\it $dim\, N>4$ and the Whitehead group $Wh(N)$ of $N$ vanishes. }
\end{enumerate}
\noindent {\bf Proof.} The proof of the existence part is the same as the proof of
the Main Theorem. Moreover it also follows from the proof of the Main Theorem that the smooth
structures $\cS$, $\cS'$ on $Q-\{q\}$ are diffeomorphic. Here
$\cS$ is the given smooth structure on $Q-\{q\}$ and  $\cS'$
is the (restriction of) normal smooth structure on
$Q-\{q\}$ with respect to $K$. \\

It remains to prove that the normal smooth structures are correct with respect to $N$. This is certainly true if $n\leq 3$ (because
of the uniqueness of differentiable structures, up to diffeomorphisms,
in dimensions $\leq 3$). We assume $n\geq 4$.\\

From the proof of the Main Theorem we see that $S=f(\sL(p,K))$ is a smooth submanifold of 
$(Q-\{q\},\cS')$ (see Lemma 1.2). Let  $\cS'_S$ be 
the smooth structure on $S$ induced by $\cS'$. Note that  
the link smoothing $h_p:(N,{\tilde{\cS}}_N)\ra(S,\cS'_S)$ 
is a diffeomorphism (see (2) above). Note also that
$\cS'_S$
is a normal smooth structure on $S$ with respect to $L$, as can be verified from the
proof of the Main Theorem applied to $S$ and $L$ (actually the link smoothings are the same; for more details see 1.3.5 in \cite{O1}).
If $n=4$ then 
Applying A.2.1 we get that $(S,\cS'_S)$ is $PL$ homeomorphic to $L$.
On the other hand condition (v) above imply that $(N,\cS_N)$ is also $PL$ homeomorphic
to $L$.
It follows that all three $(S,\cS'_S)$, $(N,{\tilde{\cS}}_N)$, $(N,\cS_N)$
are $PL$ equivalent. Since we are in dimension four 
(see \cite{Cerf}, \cite{HirschLectures}) all three spaces above
are diffeomorphic and (a) follows. \\

We prove (b). By hypothesis there is an embedding $(N,\cS_N)\ra (Q-\{q\},\cS)$
whose image lies near $q$. Since $\cS$ is diffeomorphic to $\cS'$, we get that
there is an embedding $(N,\cS_N)\ra (Q-\{q\},\cS')$ whose image lies near $q$.
But near $q$ the smooth structure $\cS'$ is diffeomorphic to the product $I_\epsilon\times S$ with smooth structure $\cS_{I_\epsilon}\times\cS'_N$, where $I_\epsilon=(0,\epsilon)$. And this product is diffeomorphic to 
$I_\epsilon\times N$ with smooth structure $\cS_{I_\epsilon}\times {\tilde{\cS}}_N$. 
Therefore there is an embedding $(N,\cS_N)\ra (I_\epsilon\times N,\cS_{I_\epsilon}\times {\tilde{\cS}}_N)$. Hence, by a standard topological
argument, 
$(N,{\tilde{\cS}}_N)$ and $(N,\cS_N)$ are smoothly h-cobordant. 
Consequently the proof of (b) is obtained by applying the h-cobordism Theorem. This
proves the Theorem.
\vspace{.3in}

\noindent {\bf \large  Appendix. Proof of Lemma 1.3 and  1.4 (1).}

In the first part of this appendix we prove Lemma 1.3,  and in the final part
we prove statement (1) of Proposition 1.4 (this is condition
{\bf C}$(4)$).
Recall that we are using the following definition: a non-degenerate $PD$ map $K\ra M$, $K$ a complex and $M$ a smooth manifold, is such that its restriction to each simplex (or cube) is smooth and
its derivative at every point is injective (\cite{MunkresLectures} p. 75).\\

Here are some ideas of the proof of Lemma 1.3.
 Recall that $f\big(\sL(q^{n-5},K)\big)$
 with the $PL$ structure induced by $K$ is a $PL$ 4-sphere
 (because $K$ induces a $PL$ triangulation on $M$).
 Also recall that the $PL$ structure on $f\big(\sL(q^{n-5},K)\big)$ induced by $\cS_5$ may not coincide with the $PL$
 structure induced by $K$. To prove Lemma 1.3 we will prove that
 these two $PL$ structures are $PL$ equivalent. The way
 to prove this is similar to the proof of condition {\bf C}$(k)$ that 
 says that all structures $\cS_k$ are
 diffeomorphic to the original $\cS$. The problem is
 that these diffeomorphisms do not respect the
 links, hence the cannot be restricted to give 
 diffeomorphisms (hence $PL$ homeomorphisms) of links. 
 But we will show that these diffeomorphisms do preserve a
 ``cube" version of the links. We introduce this concept
 in the next subsection.
 \vspace{.1in}

\noindent {\bf A.1. Cubic Links.}\\ 
Let $\epsilon >0$ and write 
$$\begin{array}{lll}I^k(\epsilon)=[0,\epsilon]^k\sbs\R^k&&
S^k(\epsilon)=\bS^{k-1}(\epsilon)\cap I^k(\epsilon)\\
C^k(\epsilon)=\D^{k}(\epsilon)\cap I^k(\epsilon)&&
F^k(\epsilon)=\{(x_1,...x_{k+1})\in I^k(\epsilon)\, ,\, x_i=1\, \,\,{\mbox{for some } i}\} \end{array}$$
\noindent where $\bS^{k-1}(\epsilon)$ and $\D^k(\epsilon)$ are the $(k-1)$-sphere and the $k$-disc of radius $\epsilon$ respectively. 
Then $S^k(\epsilon)$ is a canonical all-right spherical
simplex of dimension $k$ and radius $\epsilon$ and $I^k(\epsilon)$ is the canonical $k$-cube of side-length $\epsilon$.
Note that for $\epsilon\leq \epsilon' << \epsilon''$ we have  $C^k(\epsilon)\sbs I^k(\epsilon')\sbs C^k(\epsilon'')$.
In this situation by deforming along rays emanating from the origin we can construct a radial $PD$ self-homeomorphism
$\Theta=\Theta(\epsilon, \epsilon' , \epsilon'')$ on $C^k(\epsilon'')$ such that\, {\bf (1)} it is the identity on 
$C^k(\epsilon/2)$, and near $S^k(\epsilon'')$ \, {\bf (2)} sends $I^k(\epsilon')$ to $C^k(\epsilon)$. 
The inverse of $\Theta$ is also $PD$. \\

Let $K^k$ be a cube complex, $\sigma^i\in K$ and $x\in\dsigma^i$. Write $j=k-i$.
Recall that the $\epsilon$-link $\sL_\epsilon(\dsigma^i,K)$
and $\epsilon$-star  $\rC\sL_\epsilon(\dsigma,K)$ of $\sigma^i$ at $x$ are built by gluing copies of
$S^j(\epsilon)$ and $C^j(\epsilon)$ respectively. The $\epsilon$-{\it cubic link}\,\, $\Box\sL_\epsilon(\dsigma,K)$ of
$\sigma^i$ at $x$ is obtained by gluing copies of $F^j(\epsilon)$ in the same way as the $S^j(\epsilon)$.
Analogously, the $\epsilon$-{\it cubic cone link}\,\, $\Box\rC\sL_\epsilon(\dsigma,K)$ of
$\sigma^i$ at $x$ is obtained by gluing copies of $I^j(\epsilon)$ in the same way as the $C^j(\epsilon)$.
Then we can write   $$\rC\sL_\epsilon(\dsigma,K)\,\,\sbs\,\, \Box\rC\sL_{\epsilon'}(\dsigma,K)\,\,\sbs\,\,\rC\sL_{\epsilon''}(\dsigma,K)$$
\noindent for $\epsilon\leq\epsilon' << \epsilon''$. By gluing the maps $\Theta$ mentioned above simplexwise we
obtain a radial $PD$ self-homeomorphism $\Theta$ (we use the same letter) on $\rC\sL_{\epsilon''}(\dsigma,K)$ 
with similar properties: {\bf (1)} it is the identity on 
$\rC\sL_{\epsilon/2}(\dsigma,K)$, and  near $\rC\sL_{\epsilon''}(\dsigma,K)$ 
\, {\bf (2)} sends $\Box\rC\sL _{\epsilon'}(\dsigma,K)$ to $\rC\sL_\epsilon(\dsigma,K)$. 
Again the inverse of $\Theta$ is also $PD$.\vspace{.1in}

\noindent {\bf A.1.1. The spherical case.}
Consider the upper-half space $\{x_{n+1}>0\}\sbs \R^{n+1}$.
For each $i=1,...,n$ define the smooth function
$\theta\0{i}:\{x_{n+1}>0\}\ra  (-\pi/2,\pi/2)$ by
$\theta\0{i}(x_1,...,x_{n+1})=\tan^{-1}(\frac{x_i}{x_{n+1}})$.
Note that $\theta_i$ is the composition of the projection $(x_1,...,x_{n+1})\mapsto (x_{n+1}, x_i)$, with angle function of polar
coordinates.  Write $\theta=(\theta\0{1},...,\theta\0{n})$.
We denote the restriction of $\theta$ to the upper hemisphere 
$\bS^k_+=\bS^k\cap\{x_{n+1}>0\}$ by the same letter $\theta$. Let $p=e_{n+1}\in\bS^k_+$ (the north pole).
Then $\theta(p)=0\in\R^k$. Note that the derivative $D\theta_p$
is the identity matrix. (To see this
let $\alpha\0{1}(s)=(\sin(s),0,...,0,\cos(s))$. Then $\theta\0{1}(\alpha(s))=s$, hence $D\theta_p.e_1=(\theta\circ \alpha\0{1})'(0)=e_1$. Similarly
for $i>1$.) It follows that $\theta$ is a 
diffeomorphism near $p$ and it has an inverse $\vartheta:\B_\mu(\R^k)\ra
\bS^k$, where $\B_\mu(\R^k)$ is the ball of (small) radius
$\mu$ in $\R^k$ centered at the origin.

If $K$ is an all-right-spherical complex instead of a cube one, similar concepts
can be defined as above, with
the $\epsilon$ objects replaced by their images in $\bS^k$ by $\vartheta$,
were we take $\epsilon$ small enough.
\vspace{.1in}

\noindent {\bf A.2. Proof of Lemma 1.3.}
We write $L=\sL(q^j,K)$, $j=n-5$ and $N=f(L)$. 
Then $L$ is an all-right spherical complex $PL$ homeomorphic to the 4-sphere $\bS^4$.
From the proof of the Main Theorem it can be checked that
the smooth structure $\cS_N'$ on $N=f(L)$ given by Lemma 1.2 has an atlas of the form
$$\cA\, =\,\Big\{ \, \big(\,  h^\bullet_{\sigma^i}  \, ,\, \D^{4-i}\times\dsigma^i \,\big) \,\Big\}_{\sigma^i\in L}$$
\noindent That is, $\cA$ is a normal smooth atlas, hence
$\cS_N'$ is a normal smooth structure. 
Since $L$ is $PL$ homeomorphic to $\cS^4$, Lemma 1.3 is implied
by the following result and the fact that $PL$ homeomorphic
smooth 4-manifolds are diffeomorphic (\cite{Cerf}, \cite{HirschLectures}).\\

\noindent {\bf Proposition A.2.1.}  {\it Let $\cS'$ is a normal smooth structure
on the  manifold $N^4$ with respect to $L$, where $L$ is a 
cubulation or all-right spherical triangulation of $N$. Then $L$ is $PL$ homeomorphic
to $(N^4,\cS')$}.\\

We will give a proof of the proposition for the case of $L$
being a cube complex. The proof for the all-right spherical
case is obtained replacing cubic links by the objects in
A.1.1.\\

\noindent {\bf Proof of A.2.1.}
Since in what follows the map $f$ is not essential, 
from now on we identify $L$ and $N$ via $f$ and write $L=N$.
The identity map $L\ra N$ is not $PD$ (see remark 7 after the Main Theorem).
To prove that $(N,\cS_N')$ is $PL$ homeomorphic to $L$ we will modify the identity
$L\ra N$ to a non-degenerate $PD$ homeomorphism $\phi:L\ra (N,\cS_N')$, which implies that
$\phi:L'\ra (N,\cS_N')$ is a smooth triangulation of $(N,\cS_N')$.
Here $L'$ is a subdivision of $L$.
 We will need the following Lemma.\vspace{.1in}

\noindent {\bf Lemma A.2.2.} {\it  Let $\varphi:J\ra\bS^k$ be a smooth triangulation of
$\bS^k$ and $\psi:\bS^k\ra\bS^k$ homeomorphism which is non-degenerate $PD$ with respect to
$J$. Then $\varphi$ extends to a smooth triangulation $\varphi':J'\ra\D^{k+1}$ and $\psi$
extends to a homeomorphism $\psi':\D^{k+1}\ra\D^{k+1}$ which is 
non-degenerate $PD$ with respect to $J'$.}\vspace{.1in}

\noindent {\bf Addendum to Lemma A.2.2.} {\it We can assume that $J'$ and $\psi'$
are radial outside a small ball centered at the origin}.\vspace{.1in}

\noindent {\bf Remark.} Simply coning $\varphi$ does not work because differentiability would be lost at the origin. The idea of the
proof is to first extend $\varphi$ radially (but not to the origin).
Then, near the origin approximate this extension by a simplicial
map, and this simplicial map can be extended to the origin
using coning.\vspace{.1in}

\noindent {\bf Proof.}  Take
$K=J\times [\eta,1]$ and $f(x,t)=t\varphi(x)$, which is a non-degenerate
$PD$ homeomorphism onto its image (hence a smooth triangulation of its image \cite{MunkresLectures}, p.77). Here $\eta>0$ is small.
Let $\epsilon$ sufficiently small so that any $\epsilon$-approximation
of $f$ is also a $PD$ embedding (\cite{MunkresLectures}, Th. 8.8).
Take $K_1=J\times [\eta,2\eta]$ and $g: K'_1\ra \R^{k+1}$ the secant $\delta$-approximation of $f|_{K'_1}$ (\cite{MunkresLectures}, p. 87) where
$K_1'$ is a subdivision of $K_1$ so that $g$ is also an embedding
(\cite{MunkresLectures}, Th. 8.8).
Taking $\delta$ even smaller if necessary and using 9.8 of
\cite{MunkresLectures} (with $\epsilon$ as above) we get a subdivision $K'$ of
$J\times [\eta,1]$ and a non-degenerate $PD$ homeomorphism $h$ on $K'$ which coincides
with $f$ outside (say) $J\times [0,3\eta]$ and is simplicial on 
$K'_1$ (which is a subdivision of $J\times [\eta,2\eta]$). Let $J''$ be the subdivision of  $J\times\{\eta\}\sbs K$
induced by $K'$.
Take $J'=(\rC J'') \,\coprod _{J''}K'$ and extend $h$ simplicially on
$\rC J''\sbs J'$ to obtain $h':J'\ra \D^{k+1}$. Then
$h'$ is a non-degenerate $PD$ map and if $\epsilon$ is sufficiently small
$h'$ is also a homeomorphism. It can be checked that $J'$ depends only on the 
bounds of the derivatives of $\varphi$. Therefore we can assume that when we
apply the same argument to $\psi\circ\varphi:J\ra \bS^k$ we obtain a similar
map $h''$ defined also on the same $J'$, that is $h'':J'\ra\D^{k+1}$ is non-degenerate,
is radial outside a small neighborhood of the cone point of $J'$ and coincides
with $\psi\circ \varphi$ on $J\sbs J'$.
Finally take $\psi'=h''\circ (h')^{-1}$. This proves Lemma
A.2.2.\\

Note that, by construction, the complex $J'$ and the map $\psi'$
are radial outside a ball centered at the origin; this proves
the addendum to Lemma A.2.2.\\

We now continue the proof of Proposition A.2.1.
To simplify we write $N=(N,\cS_N')$.
Now, the identity map $L\ra N$ is already $PD$ outside the two skeleton $L^2$ of
$L$ (it is $PD$ on the interior of the 3-simplices). For $\sigma\in L$ and $\delta>0$ denote by $\sigma\0{\delta}$ the set of points
in $\sigma$ that lie at a distance $>\delta$ from $\p\sigma$. Hence 
$\sigma\0{\delta}\sbs\dsigma$. \\

Using the charts
 $\big(\,  h^\bullet_{\sigma^i}  \, ,\, \D^{4-i}\times\dsigma \,\big)$ we identify
$ \D^{4-i}\times\dsigma $ with its image by $ h^\bullet_{\sigma^i}$.
In particular, for $x\in\dsigma$, we are identifying $\D^{n-i}\times\{x\}$ 
with the $\rC\sL(\sigma,K)=\rC\sL_\epsilon(\dsigma,K)$
based at $x$, for some $\epsilon>0$.  \\

We choose $\delta>0$ small and $\epsilon >0$ even smaller to make sure that
the open sets $\D^{2}\times\sigma\0{\delta}$, for all 2-simplices $\sigma^2\in L$, have
disjoint closures.
Let $\sigma^2\in L$ be a 2-simplex in $L$ and $x\in\sigma\0{\delta}^2$. 
Now, the identity $\D^2\times\{x\}=\rC\sL (\dsigma^2,L)\ra\D^2\times\{x\}$ 
is not, in general, a $PD$ map at the origin $(0,x)\in \D^2\times\{x\}$ 
(in fact it looks like the spherical version of the example given in the
introduction; see also remark 7 in the introduction).
But using 
Lemma A.2.2
we can modify the identity $\D^2\times\{x\}=\rC\sL (\dsigma^2,L)\ra\D^2\times\{x\}$ 
near the origin $(0,x)\in \D^2\times\{x\}$ to make it
a non-degenerate $PD$ homeomorphism. Crossing with the identity $1_\dsigma$ we obtain a non-degenerate
$PD$ self-homeomorphism of $\D^2\times\sigma\0{\delta}^2$,
which is the identity outside a small neighborhood of $\sigma^2\0{\delta}$. Patching all these maps for all $\sigma^2\in L$
together with the identity we
obtain a non-degenerate $PD$ homeomorphism $\psi$  on $L(2)=(L-L^2)\cup\big(\bigcup _{\sigma^2\in L}\D^2\times\sigma^2\0{\delta}\big)$.
Note that $\psi$ is still a product map on each $\D^2\times\sigma^2\0{\delta}$. \\

Now let $\sigma^1\in L$ be a 1-simplex and $x\in\sigma^1\0{\delta'}$, where $\delta'$ is small but with
 $\delta<<\delta'$ (we may have to take $\delta$ and $\epsilon$ even smaller). Also let $\epsilon'$ be such that $\epsilon <<\epsilon'$ so that 
the $\epsilon'$-link at  $x$ is such that $\sL_{\epsilon'}(\dsigma^1,L)\sbs L(2)$, hence we can apply $\psi$ to   $\sL_{\epsilon'}(\dsigma^1,L)$.
We would like to use Lemma A.2.2 to extend $\psi$ to
a $PD$ map on $\rC\sL_{\epsilon'}(\dsigma^1,L)$ (thus also 
near   $\sigma^1\0{\delta}$) but the problem is that
$\psi$ does not map $\sL_{\epsilon'}(\dsigma^1,L)$ to itself. To correct this we work with cubic links.
Note that $\Box\sL_{\epsilon'}(\dsigma^1,L)\sbs L(2)$. And, since 
$\psi$ is still a product map on each $\D^2\times\sigma^2\0{\delta}$ we have that
now $\psi$ maps $\Box\sL_{\epsilon'}(\dsigma^1,L)$ to itself. We can now apply Lemma A.2.2 to the restriction of
$\Theta\,\psi\,\Theta^{-1}$ to $\sL_{\epsilon'}(\dsigma^1,L)$ to obtain a 
non-degenerate $PD$ extension $\psi'$ defined on $\rC\sL_{\epsilon'}(\dsigma^1,L)$. Take now
$\Theta^{-1}\, \psi'\,\Theta$ and cross this map with $1_{\sigma^1\0{\delta'}}$ to extend $\psi$ near $\sigma^1\0{\delta'}$.
Doing this for all $\sigma^1\in L$ and patching these maps with the previous $\psi$  we
obtain a non-degenerate $PD$ homeomorphism $\psi$  on $L(1)=L(2)\cup\big(\bigcup _{\sigma^1\in L}\D^3\times\sigma^1\0{\delta}\big)$.
Note that $\psi$ is still a product map on each $\D^2\times\sigma^2\0{\delta''}$ and on each $\D^3\times\sigma^1\0{\delta'}$ (where $\delta''$ is slightly smaller than
$\delta$). Therefore $\psi$ maps $\Box\sL_{\epsilon'''}(\dsigma^0,L)$ to itself, for a suitable choice of $\epsilon'''$.
To extend $\psi$ to the whole of $L$ we proceed one step further (now for 0-simplices $\sigma^0$) in a similar way in the case for 1-simplices
$\sigma^1$.
This proves Proposition A.2.1 and Lemma 1.3.\\

The argument used in the proof above gives the following corollaries. 
In the first corollary we use the notation in the Main Theorem and its proof.\vspace{.1in}

\noindent {\bf Corollary A.2.3.} {\it Let $f:K\ra (M^4,\cS)$ be a PD non-degenerate
homeomorphism.
Let $\cS'=\cS_5$ (as in the proof of the Main Theorem). Then $(M,\cS)$ and $(M,\cS')$
are diffeomorphic.}\vspace{.1in}

\noindent {\bf Proof.} We have that $\cS'$ is a normal smooth structure
on $M$, with respect to $K$. By Proposition A.2.1 $(M,\cS')$ is $PL$ homeomorphic
to $K$. On the other hand, since $K$ is a smooth cubulation 
(or all-right spherical triangulation) of 
$(M,\cS)$, we have that $K$ is $PL$ homeomorphic to
$(M,\cS)$. Therefore   $(M,\cS)$ and $(M,\cS')$ are
$PL$ homeomorphic, and, since we are in dimension 4 (\cite{Cerf},
\cite{HirschLectures}), it follows that
$(M,\cS)$ and $(M,\cS')$ are
diffeomorphic. This proves the corollary.\\

We can also generalize the arguments in the proof of Lemma 1.3 to higher dimensions. 
Suppose $K$ is a (cube of all-right spherical) complex of dimension $n$. We can choose the $\epsilon$'s and 
$\delta$'s above properly and define 
$K(n-2)=(K-K_{n-2})\cup\big(\bigcup _{\sigma^{n-2}\in K}\D^2\times\sigma^{n-2}\0{\delta}\big)$ and
$K(k)=K(k+1)\cup\big(\bigcup _{\sigma^j\in L}\D^{n-j}\times\sigma^j\0{\delta}\big)$. From the definitions we have the equalities
{\small$$
\begin{array}{rcccl} K-K(n-2)&=&K_{n-2}-\bigcup\0{\sigma^{n-2}\in K}\sigma\0{\delta}^{n-2}&\supset&K_{n-3}\\\\
K-K(k)&=&\big(K-K(k+1)\big)-\bigcup\0{\sigma^{k}\in K}\sigma\0{\delta}^{k}&\supset&K_{k-1}\end{array}
$$}

\noindent  where the inclusion on top holds because
the $\sigma^{n-2}\0{\delta}$ are disjoint from $K_{n-3}$,
and the inclusion on the second row follows from the first and induction. 
Note that these inclusions are homotopy equivalences. Hence we get\\

\noindent {\bf (A.2.4)}\hspace{2in} $K(k+1)\, \sbs \,K-K_k$\\

\noindent Moreover, it is straightforward to verify that this
inclusion is a homotopy equivalence.\\

\noindent {\bf Corollary A.2.5.} {\it Let $f:K\ra (M^n,\cS)$ be a PD non-degenerate
homeomorphism.
Let $\cS_k$ be as in the proof of of the Main Theorem satisfying condition
} {\bf S($k$)}. {\it Then there is a non-degenerate $PD$ homeomorphism
$\psi:K(k+1)\ra \big(f(K(k+1)),\cS_k\big)$. In particular we have that (take $k=-1$) $(K,\psi)$ is a smooth cubulation of $(M,\cS')$, for some $\psi$.}\\

\noindent {\bf Proof.} Take $K=L$ in the proof of Lemma 1.3 above, but
now assume $K$ of arbitrary dimension. 
Using exactly the same arguments (i.e. use Lemma A.2.2 and cubic links) given there we can extend the induction to all values of $k$. This proves the corollary.\vspace{.1in}

We now prove statement (1) of Proposition 1.4.\\

\noindent {\bf Proof of C($4$).}
The smooth structure $\cS_4$ is defined on $M-M_{n-4}$ and we have to prove
that it extends to a smooth structure $\cS_4'$ on $M$, and $\cS_4'$ is
diffeomorphic to $\cS$.\vspace{.1in}

Write $A=M-M_{n-4}$.
Let $L(j)$ be as in corollary A.2.5 and write $B=L(n-3)\sbs K$, $D=f(B)\sbs M$. 
Note that $D\sbs A$. Since the inclusion in A.2.4 is a
homotopy equivalence we have:\\

\noindent {\bf (A.2.6)}\hspace{.2in} {\it the inclusion $D\hookrightarrow A$ is a
homotopy equivalence.}\\

By corollary A.2.5 there is a non-degenerate $PD$ homeomorphism
$\psi:B\ra (D,\cS_4|_{D})$. We have that $\psi$ extends to a homeomorphism
(we use the same letter) $\psi:K\ra M$.\\

\noindent {\bf Remark.}
The last statement follows from the inductive proof of Corollary A.2.5, taking $j=-1$. That is, the $\psi$ constructed for $K(j)$ extends the $\psi$ constructed in the previous step for $K(j+1)$,
but on a slightly smaller smaller $K(j)$; see the proof of Lemma 1.3.
But we can also construct this extension using Alexander's trick (beyond codimension 3) instead of Lemma A.2.2.\vspace{.1in}

Write $f'=f|_B$ and note that, by hypothesis, 
$f':B\ra (D,\cS)$ is also a non-degenerate $PD$ homeomorphism.
Consider the smooth structures $\cS_\psi=\psi^*\cS_4$ and
$\cS_f=(f')^*\cS=(f^*\cS)|_B$ on $B$. The $PL$ structure on $B$
(induced by $K$) is thus Whitehead compatible with both smooth
structures $\cS_\psi$ and $\cS_f$. But $A=M-M_{n-4}$ has the homotopy type
of a 3-complex hence, by A.2.5, so do $D$ and $B$. Since $PL/O$ is 6-connected
it follows from the theory of smoothings of $PL$-manifolds (see Theorem
4.2 in the second essay in \cite{HirschLectures}) that $\cS_\psi$ and $\cS_f$ are concordant.
Consequently, pushing forward everything to $D$ by $\psi$ we have that the smooth structures $\cS_4|_D$ and $\cS''=\psi_*\cS_f$ are concordant.   Note that
$\cS''=\Big((f\circ \psi^{-1})^*\cS\Big)_{D}$ hence we get that\\

\noindent {\bf (A.2.7)}\,\,\,\,\,\,\,
{\it  the smooth structures $\cS_4|_D$ and $\Big((f\circ \psi^{-1})^*\cS\Big)_{D}$ on $D$ are concordant}\\

\noindent It follows from A.2.6, A.2.7 and
the theory of smoothings of topological manifolds (see the Classification Theorem 10.1,
p .194, in \cite{KiSi} and its naturality for restrictions) that the
smooth structures $\cS_4$ and $\Big((f\circ \psi^{-1})^*\cS\Big)_{A}$ on $A$
are concordant. Therefore we can find a self-homeomorphism $g$ on $A$
such that $\Big((f\circ \psi^{-1})^*\cS\Big)_{A}=g^*\cS_4$.
Moreover we can assume (see Theorem 4.1, p. 25, in \cite{KiSi}) that
that $d\0{M}(x,g(x))\ra 0$, as $x\ra\p A$ (here $d\0{M}$ is any metric on $M$
inducing the topology on $M$). Therefore we can extend $g$ to a self-homeomorphism
on $M$ by defining $g(x)=x$, for $x\notin A$. To finish the proof just take
$\cS'_4=(f\circ \psi^{-1}\circ g^{-1})^*\cS$, which
is defined on the whole $M$ and extends $\cS_4$. This proves statement {\bf C}(4).\\

Pedro Ontaneda

SUNY, Binghamton, N.Y., 13902, U.S.A.

\end{document}